\documentclass[12pt,twoside]{article}
\usepackage{amsmath,amssymb,amsthm}
\usepackage{graphicx}

\usepackage{datetime}
\date{}

\newtheorem{theorem}{Theorem}[section]
\newtheorem{proposition}[theorem]{Proposition}
\newtheorem{lemma}[theorem]{Lemma}
\newtheorem{corollary}[theorem]{Corollary}
\newtheorem{example}[theorem]{Example}

\newtheorem{problem}[theorem]{\bf Problem}

\title{\textbf{Structure and Applications of $(p, q)$-Derivations in Algebra}}

\author{
\textbf{Mohammad Javad Mehdipour}\thanks{Corresponding author}\\
Department of Mathematics\\
Shiraz University of Technology\\
Shiraz 71555-313, Iran\\
\texttt{mehdipour@sutech.ac.ir}
\and
\textbf{Narjes Salkhordeh}\\
Department of Mathematics\\
Shiraz University of Technology\\
Shiraz 71555-313, Iran\\
\texttt{n.salkhordeh@sutech.ac.ir}
}

\begin{document}

\maketitle

\begin{abstract}
In this paper, we investigate $(p, q)$-derivations in both general algebras and Banach algebras. Our main result extends the Singer-Wermer theorem to $(p, q)$-derivations, proving that their ranges are contained in the radical of the algebra. As a direct consequence, we prove that the range of any left derivation lies within the radical of the algebra even without assuming continuity. This improves a result due to Bra\v{s}er and Vukman. For Banach algebras, we show that primitive ideals remain invariant under bounded $(p, q)$-derivations. Finally, we study $(p, q)$-derivations of group algebras and give an answer to the question posed.
\end{abstract}



\section{Introduction}

Derivations and their variants play a central role in understanding the structure of rings and algebras. These concepts have been studied extensively, with particular attention given to Jordan derivations, left and right derivations, and other functional identities; see, e.g., \cite{bm, bv1, her, hv, 13}. These investigations have established deep connections between derivations and various structural properties, such as primeness, semiprimeness, torsion conditions, and the behavior of radicals and ideals.

In this context, Vukman and Kosi-Ulbl \cite{10} established a fundamental result. They proved that if $R$ is a semiprime ring with identity under appropriate characteristic hypotheses, and $d: R\rightarrow R$ is an additive map satisfying
$$
d(x^{p+q+1})=(p+q+1)\,x^{p}d(x)\,x^{q}\qquad(x\in R),
$$
for distinct non-negative integers $p,q$ with $p+q\neq0$, then $d$ is a derivation whose image lies in the center $Z(R)$. In the course of their argument they encountered additive maps satisfying
\begin{eqnarray}\label{*}
(p+q)d(x^2)=2pd(x)x + 2qxd(x)\qquad(x\in R).
\end{eqnarray}
Clearly, when $q=0$, identity (\ref{*}) reduces to the Jordan left condition. By \cite{vukman},  Jordan left derivations on 2- and 3-torsion-free semiprime rings are derivations whose images lie in the center. These observations led to the conjecture that any additive map satisfying (\ref{*}) must itself be a derivation mapping $R$ into $Z(R)$. To solve this conjecture, they considered the more general weighted identity
\begin{eqnarray}\label{**}
(p+q)d(xy)=2pd(x)y + 2qxd(y)\qquad(x,y\in R),
\end{eqnarray}
and proved that $d=0$ whenever $R$ is semiprime under the relevant characteristic hypotheses. Let us remark that the conjecture was solved by Kosi-Ulbl and Vukman $\cite{9}$. In what follows we refer to maps satisfying (\ref{**}) as $(p,q)$-derivations, and we investigate their structure when $R$ is an algebra.

Note that if $p=0$ and $q\neq0$, then a $(p,q)$-derivation is a left derivation, while if $q=0$ and $p>0$, it becomes a right derivation. Moreover, if we allow $p=q>0$, the identity reduces to that of an ordinary derivation. We denote by $\hbox{Der}_{p,q}(A)$ the set of all $(p,q)$-derivations of an algebra $A$, and by $\hbox{Der}(A)$ the set of all ordinary derivations. For a normed algebra $A$, the subsets consisting of bounded elements are denoted by $\hbox{Der}_{p,q}^c(A)$ and $\hbox{Der}^c(A)$, respectively.

The purpose of this paper is to extend the study of $(p,q)$-derivations beyond semiprime rings, focusing on normed algebras and group algebras. In Section 2, we investigate $(p,q)$-derivations of general algebras and prove that the range of any $(p,q)$-derivation of an algebra $A$ is contained in the radical $\hbox{rad}(A)$. As a consequence, we show that every left derivation maps $A$ into $\hbox{rad}(A)$, which demonstrates that the continuity assumption in Theorem 2.1 of \cite{bv1} is unnecessary. Furthermore, we prove that the primitive ideals of a Banach algebra remain invariant under its bounded $(p,q)$-derivations. In Section 3, we investigate $(p, q)$-Jordan derivations of algebras associated with locally compact groups and provide an answer to the question posed in \cite{am}.

\section{\normalsize\bf Singer-Wermer and related results}

We commence this section with a lemma.

\begin{lemma}\label{1} Let $ I$ be an ideal of  an algebra $A$ and $d\in\hbox{Der}_{p, q}(A)$. Then the following statements hold.

\emph{(i)} For every $ a_1, a_2\in A$ and $ n \in {\Bbb N} $
$$ d^{n} (a_1a_2) = \frac{2^{n}}{(p+q)^{n}} \sum \limits_{k=0}^n \displaystyle{n\choose k} p^{n-k} q^{k} d^{n-k} (a_1) d^{k} (a_2),$$
with the convention $0^0=1$ and $d^{0}(a)= a $ for all $a \in A$.

\emph{(ii)} For every $\iota\in I $ and $ n \in {\Bbb N}$,
$$ d^{n} (\iota^{n}) - n! \left( \frac{2}{p+q} \right) ^{\frac{n^{2}+n-2}{2}} q^{\frac{n^{2}-n}{2}} p^{n-1} d(\iota)^{n} \in I. $$ 
\end{lemma}
{\it Proof.} 
(i) Assume the formula holds for some fixed \( n \in {\Bbb N} \). For any \( a_1, a_2 \in A \)
$$
d^{n+1} (a_1a_2) = d ( d^{n} (a_1a_2) ) =d(\frac{2^{n}}{(p+q)^{n}} \sum \limits_{k=0}^n \displaystyle{n\choose k} p^{n-k} q^{k} d^{n-k} (a_1) d^{k} (a_2)).
$$
Applying the linearity of \( d \) and the \((p,q)\)-derivation identity to each product $d^{n-k} (a_1) d^{k} (a_2)$, we obtain
\begin{eqnarray} 
d^{n+1} (a_1a_2) &=& \frac{2^{n+1}}{(p+q)^{n+1}} ( \sum \limits_{k=0}^n \displaystyle{n\choose k} p^{n+1-k} q^{k}d^{n+1-k} (a_1) d^{k} (a_2)\nonumber \\
&+& \sum \limits_{k=1}^{n+1} \displaystyle{n\choose k-1} p^{n+1-k} q^{k}d^{n+1-k} (a_1) d^{k} (a_2)). \nonumber 
\end{eqnarray}
Simplifying the expression,  we have
\begin{eqnarray}
d^{n+1} (a_1a_2) &=& \frac{2^{n+1}}{(p+q)^{n+1}} ( \displaystyle{n\choose 0} p^{n+1} d^{n+1} (a_1) a_2 \nonumber\\ &+& \sum \limits_{k=1}^n (\displaystyle{n\choose k}+\displaystyle{n\choose k-1}) p^{n+1-k} q^{k}d^{n+1-k} (a_1) d^{k}(a_2) \nonumber \\
&+& \displaystyle{n\choose n} q^{n+1} a_1 d^{n+1} (a_2) ). \nonumber \end{eqnarray}
It follows that
\begin{eqnarray}
d^{n+1} (a_1a_2) &=& \frac{2^{n+1}}{(p+q)^{n+1}}( \displaystyle{n+1\choose 0} p^{n+1} d^{n+1} (a_1) a_2 \nonumber \\
&+& \sum \limits_{k=1}^n \displaystyle{n+1\choose k} p^{n+1-k} q^{k}d^{n+1-k} (a_1) d^{k}(a_2) \nonumber \\
&+& \displaystyle{n+1\choose n+1} q^{n+1} a_1 d^{n+1} (a_2) ) \nonumber \\
&=& \frac{2^{n+1}}{(p+q)^{n+1}} (\sum \limits_{k=0}^{n+1} \displaystyle{n+1\choose k} p^{n+1-k} q^{k} d^{n+1-k} (a_1) d^{k} (a_2)), \nonumber 
\end{eqnarray}
which proves the formula for $n+1$. Thus (i) holds by induction.
%

(ii) Assume the claim holds for some fixed $ n > 2 $. For every $\iota\in I$, we have
 $$
d^{n+1} (\iota^{n+1})+I = \frac{1}{p+q} (p+q) d ( d^{n} (\iota^{n+1}) ) +I.
$$
Using part (i) to expand $d ( d^{n} (\iota^{n+1}) )$, we obtain
\begin{eqnarray}
d^{n+1} (\iota^{n+1})+I &=& \frac{2^{n+1}}{(p+q)^{n+1}} ( \sum \limits_{k=0}^n \displaystyle{n\choose k} p^{n+1-k} q^{k}d^{n+1-k} (\iota) d^{k} (\iota^{n}) \nonumber \\
&+& \sum \limits_{k=0}^n \displaystyle{n\choose k} p^{n-k} q^{k+1}d^{n-k} (\iota) d^{k+1} (\iota^{n}) ) +I\nonumber \\
&=& \frac{2^{n+1}}{(p+q)^{n+1}} ( \displaystyle{n\choose n} p q^{n}d (\iota) d^{n} (\iota^{n}) + \displaystyle{n\choose n-1} p q^{n}d (\iota) d^{n} (\iota^{n}) + I ). \nonumber 
\end{eqnarray}
Therefore, the expression above simplifies to the following.
 \begin{eqnarray}
d^{n+1} (\iota^{n+1})+I &=& \frac{2^{n+1}}{(p+q)^{n+1}} (n+1) p q^{n}d (\iota) d^{n} (\iota^{n}) +I \nonumber \\ 
&=& (n+1)! ( \frac{2}{p+q})^{\frac{(n+1)^{2}+(n+1)-2}{2}} q^{\frac{(n+1)^{2}-(n+1)}{2}} p^{n} d(\iota)^{n+1} + I, \nonumber 
\end{eqnarray}
which is exactly the desired identity for $n+1$. This completes the induction and the proof of (ii).$\hfill\square$\\

Singer-Wermer \cite{14} proved that the range of bounded derivations of a commutative Banach algebra is a subset of its radical. Thomas \cite{15} showed that this result even holds for unbounded derivations. Some authors have studied this result for non-commutative Banach algebras \cite{bm, 11, mm}.
 We now establish an analogous result for \((p,q)\)-derivations. First, recall that $\hbox{nil}(A)$ denotes the intersection of all prime ideals of \(A\), and satisfies $\hbox{nil}(A)\subseteq \hbox{rad}(A)$.

\begin{theorem}\label{2} Every $(p,q)$-derivation of an algebra $ A $, maps the algebra into $\emph{nil}(A)$. In particular, if $ A $ is semisimple, then $\hbox{Der}_{p, q}(A)=\{0\}$.
\end{theorem} 
{\it Proof.} Let $d\in\hbox{Der}_{p, q}(A)$. Fix $ a_1, a_2, a_3 \in A$ and non-negative integers $ k, n$. By Lemma \ref{1}, we have
\begin{eqnarray*}\label{r6}
d^{n+k} (a_1a_2a_3) &=& \frac{2^{n+k}}{(p+q)^{n+k}} ( \sum \limits_{j=0}^{n+k} \displaystyle{n+k\choose j} p^{n+k-j} q^{j} d^{n+k-j} (a_1) d^{j} (a_2a_3)) \nonumber \\
&=& \frac{2^{n+k}}{(p+q)^{n+k}} ( \sum \limits_{j=0}^{k-1} \displaystyle{n+k\choose j} p^{n+k-j} q^{j} d^{n+k-j} (a_1) d^{j} (a_2a_3) \\
&+& \displaystyle{n+k\choose k} p^{n} q^{k} d^{n} (a_1) d^{k} (a_2a_3) 
+ \sum \limits_{j=k+1}^{n+k} \displaystyle{n+k\choose j} p^{n+k-j} q^{j} d^{n+k-j} (a_1) d^{j} (a_2a_3)).\nonumber 
\end{eqnarray*}
We also have
$$
d^{n}(a_1) d^{k} (a_2a_3) = \frac{2^{k}}{(p+q)^{k}} ( d^{n}(a_1) a_2 d^{k}(a_3) + d^{n} (a_1) \sum \limits_{j=0}^{k-1} \displaystyle{k\choose j} p^{k-j} q^{j} d^{k-j} (a_2) d^{j} (a_3)). 
$$
Assume now that $ P $ is a minimal prime ideal of  $A$. By the argument used in the proof of the lemma of \cite{mm}, we see that $d(P)$ is a subset of $P$. So the induced $(p, q)$-derivation $D: A/P\rightarrow A/P$ defined by
$$ D(a+P) = d(a) + P $$
is well-defined. From Theorem 6 of $\cite{10}$ we conclude that $ D = 0 $. That is, $ d(a)\in P $ for all $ a \in A $. Therefore $ d(x)\in\hbox{nil}(A)\subseteq \hbox{rad}(A) $.
$\hfill\square$\\

Bre\v{s}ar and Vukman \cite{bv1} proved that if $d$ is a  continuous left derivation of a Banach algebra $A$, then $d$ maps $A$ into $\hbox{rad}(A)$. Theorem \ref{2} shows that the continuity requirement may be removed and the result  extended to arbitrary normed algebras. 

\begin{corollary} Let $d$ be a left derivation of an algebra $A$. Then $d$ maps $A$ into $\emph{rad}(A)$.
\end{corollary}

Sinclair $\cite{13}$ proved that primitive ideals of a Banach algebra are invariant under bounded derivations. We prove this result for bounded $(p, q)$-derivations.

\begin{theorem}\label{3m}  Let $A$ be a Banach algebra and $d\in\hbox{Der}_{p, q}^c(A)$ . Then primitive ideals are invariant under $ d $.
\end{theorem}
{\it Proof.}
Let $ P $ be a primitive ideal of  $ A $, and let $p$ be a non-zero integer. Then for every $x\in P$ and $a\in A$, we have
\begin{eqnarray*}
(d(x)+P)(a+P)&=& d(x)a+P\\
&=&\frac{1}{2p}((p+q)d(xa)-2qxd(a))+P\\
&=&d(\frac{p+q}{2p} xa)+P.
\end{eqnarray*}
Set $z=\frac{p+q}{2p} xa$. Then from Lemma $\ref{1}$ we infer that
$$ d^{n} (z^{n}) + P = n! (\frac{2}{p+q}) ^{\frac{n^{2}+n-2}{2}} q^{\frac{n^{2}-n}{2}} p^{n-1} d(z)^{n} + P. $$
This shows that
\begin{eqnarray} 
\| ( d(z)+ P)^{n} \|^{\frac{1}{n}} &=& \left( \frac{(p+q)^{\frac{n^{2}+n-2}{2}}}{n! 2^{\frac{n^{2}+n-2}{2}} q^{\frac{n^{2}-n}{2}} p^{n-1}}\right)^{\frac{1}{n}} \| d^{n} (z^{n}) + P \|^{\frac{1}{n}} \nonumber \\
&\leq & \left( \frac{(p+q)^{\frac{n^{2}+n-2}{2}}}{n! 2^{\frac{n^{2}+n-2}{2}} q^{\frac{n^{2}-n}{2}} p^{n-1}}\right)^{\frac{1}{n}}\|d \| \|z\|. \nonumber
\end{eqnarray}
Thus $ r(d(z) + P) = 0 $, where $r(.)$ is the spectral radius. So 
$$
r((d(x)+P)(a+P))=0
$$
for all $x\in P$ and $a\in A$. It follows that
$
d(x)+P
$ is an element of 
$$\hbox{rad}(A/P)= P;$$
see for example Proposition 1, page 126 of \cite{bd}. Therefore, $d(x)\in P$.
$\hfill\square$\\

Let $A$ be an algebra and $a\in A$. Then $b\in A$ is called a {\it quasi-inverse} of $a$ if $$ab=ba=a+b.$$ 
The set of all quasi-invertible elements of $A$ is denoted by $\hbox{q-Inv}(A)$. 
Let us recall that a scalar $\lambda\in\Bbb{C}$ is an $\emph{eigenvalue}$ of a linear operator $T: A\rightarrow A$ if there exists a non-zero element $a\in A$ such that $T(a)=\lambda a$. In this case, $a$ is called an $\emph{eigenvector}$ for $\lambda$. We denote the set of all eigenvectors of $T$ for scalar $\lambda$ by $\hbox{E}_\lambda(T)$. Clearly, $$\hbox{E}_0(T)=\hbox{ker}(T)\setminus\{0\}.$$
It is easy to see that if $\lambda$ is non-zero, then $$\hbox{E}_\lambda(T)\subseteq\hbox{rang}(T);$$ in particular, if $T$ is a $(p,q)$-derivation, then by Theorem \ref{2}, 
$$
\hbox{E}_\lambda(T)\subseteq\hbox{rad}(A)\subseteq\hbox{q-Inv}(A);
$$
see also \cite{bd}.
For fixed $n\in{\Bbb N}$, let $Q_n(A)$ be the set of all $a\in A$ such that $a^n=0.$ Note that every element of $Q_n(A)$ is nilpotent.

Mathieu and Murphy $\cite{11}$ proved that elements of $\hbox{E}_\lambda(d)$ are nilpotent, when $d$ is a bounded derivation. We establish an analogue of this result for any $ (p,q)$-derivations.

\begin{theorem}\label{4}  
Let $A$ be an algebra, $d\in\hbox{Der}_{p, q}(A)$, and $\lambda$ be a non-zero scalar. Then the following statements hold.

\emph{(i)} $\emph{E}_\lambda(d)$ is contained in $Q_3(A)$.

\emph{(ii)} If $a\in\emph{E}_\lambda(d)$, then $-a-a^2$ is the quasi-inverse of $a$.

\emph{(iii)} If $A$ is unital or is an algebra with no non-zero nilpotents of index $\leq 3$, then $\emph{E}_\lambda(d)$ is the empty set.
\end{theorem}
{\it Proof.} (i) Let $ a \in\hbox{E}_\lambda(d)$. Then $ d(a) = \lambda a $. Thus
\begin{eqnarray}
(p+q) d (a^{2}) &=& 2p d(a) a + 2q a d(a) \nonumber \\
&=& 2(p+q)\lambda a^{2}. \nonumber
\end{eqnarray}
Hence $ d (a^{2})= 2 \lambda a^{2} $. Thus
\begin{eqnarray}\label{r8}
d(a^{3} ) &=& d(a^{2} a ) \nonumber \\
&=& \frac{2}{p+q} (p d(a^{2}) a + q a^{2} d(a))\\
&=& \frac{2}{p+q}( 2p+ q) \lambda a^{3}. \nonumber 
\end{eqnarray}
We also have
\begin{eqnarray}\label{r9}
d(a^{3} ) &=& d(aa^{2} ) \nonumber \\
&=& \frac{2}{p+q} (p d(a) a^{2} + q a d(a^{2})) \\
&=& \frac{2}{p+q}( p+ 2q) \lambda a^{3}.\nonumber 
\end{eqnarray}
In view of $(\ref{r8})$ and $(\ref{r9})$, we have
$$ ( p+ 2q) \lambda a^{3} = ( 2p+ q) \lambda a^{3}. $$
If $ \lambda \neq 0 $, then $ a^{3} = 0 $. Thus $ a\in Q_3(A)$. 

(ii) For every $a \in\hbox{E}_\lambda(d)$, we have
\begin{eqnarray*}
a(-a-a^2)&=& -a^2-a^3=-a^2\\
&=&a+(-a-a^2)\\
&=&(-a-a^2)a.
\end{eqnarray*}
 So $-a-a^2$ is the quasi-inverse of $a$. 
 
(iii) Assume now that $A$ is unital with the identity $1$. Suppose that $\hbox{E}_\lambda(d)\neq\emptyset$ and choose an element $a\in\hbox{E}_\lambda(d)$. Note that $a$ is non-zero. By (i), $a^3=0$ and so 
 \begin{eqnarray}\label{zx}
 (1-a)(1+a+a^2)=1.
 \end{eqnarray}
 Thus 
 \begin{eqnarray}\label{zx1}
 (p+q)d(1)&=&(p+q)d((1-a)(1+a+a^2))\nonumber\\
 &=&2pd(1-a)(1+a+a^2)+2q(1-a)d(1+a+a^2).
 \end{eqnarray}
 But, 
\begin{eqnarray*}
(p+q)d(a^2)&=&2pd(a)a+2qad(a)\\
&=&2p\lambda a^2+2q\lambda a^2\\
&=&2\lambda(p+q)a^2.
\end{eqnarray*}
This shows that $d(a^2)=2\lambda a^2$. One can prove that $d(1)=0$. Hence (\ref{zx1}) can be rewriten to the following.
$$
2p(-\lambda a)(1+a+a^2)+2q (1-a)(\lambda a+2\lambda a^2)=0.
$$
It follows that 
\begin{eqnarray*}
2\lambda(q-p)(a+a^2)=0.
\end{eqnarray*}
Since $\lambda$ and $p-q$ are non-zero, we have $a+a^2=0$. From this and (\ref{zx}) we infer that 
$$
(1-a)(1+0)=1.
$$
Consequently, $1-a=1$ and so $a=0$, a contradiction.

Finally, let $A$ be an algebra with no non-zero nilpotents of index $\leq 3$. If $a \in E_\lambda(d)$, then (i) shows that $a^3 = 0$, and hence $a = 0$. This contradiction implies that $\hbox{E}_\lambda(d)=\emptyset.$
$\hfill\square$\\

As an immediate consequence of Theorem \ref{4} we have the following result.

\begin{corollary}
Let $A$ be an algebra and $d\in\hbox{Der}_{p, q}(A)$. Then the following statements hold.

\emph{(i)} If $A$ is a semisimple or semiprime algebra, then $\emph{E}_\lambda(d)=\emptyset$ for all non-zero scalars $\lambda$.

\emph{(ii)} If $d$ can extend to the unitization of $A$ as a $(p, q)$-derivation, then $\emph{E}_\lambda(d)=\emptyset$ for all non-zero scalars $\lambda$. 
\end{corollary}

The following example shows that the bound $3$ in Theorem \ref{4} (i) is optimal. It also shows that there is no relation between the eigenvalues of a $(p,q)$-derivation and the values of $p$ and $q$.

\begin{example}{\rm Let $a$ be a non-zero matrix such that $a^2\neq 0$ but $a^3=0$. Let $A$ be the algebra generated by $\{a\}$. Let us remark that every element of $A$ has form of $\alpha a+\beta a^2$, where $\alpha, \beta\in {\Bbb C}$. Choose a non-zero scalar $\lambda$. Define $d: A\rightarrow A$ by 
$
d(a)=\lambda a$ and $d(a^2)=2\lambda a^2,
$
and extend linearly to all of $A$. Take arbitrary $a_1=\alpha a+\beta a^2$ and $a_2=\eta a+\delta a^2$. Their product is

$$
a_1a_2 = (\alpha a+\beta a^2)(\eta a+\delta a^2) = \alpha\eta a^2.
$$
It follows that
$$
d(a_1a_2)=d(\alpha\eta a^2)=2\alpha\eta\lambda a^2.
$$
On the other hand,
$$
d(a_1)=\alpha\lambda  a+2\beta\lambda  a^2\qquad\hbox{and}\qquad d(a_2)=\eta\lambda  a+2\delta\lambda  a^2.
$$
Assume now that $ p $ and $ q $ are distinct non-negative integers satisfying $p+q\neq 0$. Then
$$
2p\,d(a_1)a_2+2q\,a_1 d(a_2)=2(p+q)\,\alpha\eta\lambda  a^2=(p+q)d(a_1a_2).
$$
Therefore, $d$ is a $(p,q)$-derivation. By construction $d(a)=\lambda  a$. Thus $a\in\hbox{E}_\lambda (d)$. Also, $a^3=0$ while $a^2\neq0$, so the nilpotency index of $a$ is exactly $3$. Since $p$, $q$ and $\lambda$ are arbitrary, there is no relation between the eigenvalues of $d$ and the values of $p$ and $q$. 
}
\end{example}

For a Banach algebra $A$, let us recall that the first Arens product ``$\diamond$" on $A^{**}$, the second dual of $A$, is defined by 
$$
\langle F_1\diamond F_2, f\rangle=\langle F_1, F_2f\rangle,
$$
where $\langle F_2f, a\rangle=\langle F_2, fa\rangle$, in which 
$\langle fa, x\rangle=\langle f, ax\rangle$ for all $F_1, F_2\in A^{**}$, $f\in A^*$ and $a, x\in A$; see for example $\cite{6}$.

\begin{proposition}\label{6} Let $A$ be a Banach algebra. Then the following statements hold.

\emph{(i)} If $ A $ has a right identity, then $\hbox{Der}_{p, q}(A)=\{0\}$.

\emph{(ii)} If $ A $ has a bounded left approximate identity, then $\hbox{Der}_{p, q}^c(A)=\{0\}$.
\end{proposition}
{\it Proof.} (i) Let $E$ be a right identity of $ A $ and $d\in\hbox{Der}_{p, q}(A)$. Then 
\begin{eqnarray}\label{1402t}
(p+q)d(E)=(p+q)d(E\diamond E)=2p d(E)+2q E\diamond d(E).
\end{eqnarray}
Hence
$$
(p+q)E\diamond d(E)=2pE\diamond d(E)+2qE\diamond E\diamond d(E).
$$
This shows that $E\diamond d(E)=0$. From this and (\ref{1402t}) we implies that 
$ d(E) = 0$. Thus
$$ (p+q) d(x)= (p+q) d(x\diamond E)= 2p d(x) $$
for all $ x \in A $. Therefore, $ d=0 $.

(ii) Let $d\in\hbox{Der}_{p, q}^c(A)$. Then for every $ a\in A $ and $f\in A^*$, we have
\begin{eqnarray*}
(p+q) d^{*}(f) a= 2p f d(a) + 2q d^{*}(fa).
\end{eqnarray*}
So, if $F_2 \in A^{**} $, then
\begin{eqnarray*}
(p+q) F_2d^{*}(f)= 2p\;  d^{*} (F_2f) + 2q  d^{**} (F_2) f.
\end{eqnarray*}
Thus for every $ F_1 \in A^{**} $ we obtain
\begin{eqnarray*}
(p+q) d^{**} (F_1\diamond F_2)=2p d^{**} (F_1)F_2 + 2q F_1 d^{**} (F_2).
\end{eqnarray*}
That is, $ d^{**} $ is a $(p,q)$-derivation of $A^{**}$.
Since $A$ is a Banach algebra with a bounded left approximate identity, $ A^{**} $ has a right identity. By (i), $ d^{**} = 0 $ on $A^{**}$. Therefore $ d = 0 $ on $A$. 
$\hfill\square$\\

In the follows let $\frak{G}$ be a Hausdorff locally compact group, and let $L^1(\frak{G})$ and $M(\frak{G})$ denote the Banach algebras defined in \cite{8}. It is well known that $L^1(\frak{G})$ is a semisimple Banach algebra with a bounded approximate identity; hence $L^1(\frak{G})^{**}$, equipped with the first Arens product, is a Banach algebra admitting a right identity. Likewise, $M(\frak{G})$ is a semisimple unital Banach algebra.

Let $L_0^\infty(\frak{G})$ (resp. $M_*(\frak{G})$) be the subspace of $L^1(\frak{G})^*$ (resp. $M(\frak{G})^*$) consisting of all functions that vanish at infinity. Then $L_0^\infty(\frak{G})^*$ and $M_*(\frak{G})^*$ are Banach algebras under the first Arens product. One can show that $L_0^\infty(\frak{G})^*$ has a right identity, whereas $M_*(\frak{G})^*$ is unital; for further details on $L_0^\infty(\frak{G})^*$ see \cite{am1, lp,mn,mn1,ms}, and for $M_*(\frak{G})^*$ see \cite{mm1, m}.

Finally, let $A(\frak{G})$ and $B(\frak{G})$ denote the Fourier and Fourier-Stieltjes algebras of $\frak{G}$. It is well known that $A(\frak{G})$ has a bounded approximate identity and that $B(\frak{G})$ is unital. Combining these facts with Theorem \ref{2} and Proposition \ref{6}, we obtain the following result.

\begin{corollary}\label{am} The following statements hold.

\emph{(i)} If $\frak{A}$ is either $L^1(\frak{G})$ or $F(\frak{G})$, then the space $\hbox{Der}_{p, q}^c(\frak{A})$ is trivial.

\emph{(ii)} If $\frak{A}$ is one of $M(\frak{G})$, $B(\frak{G})$, $L_0^\infty(\frak{G})^*$, $L^1(\frak{G})^{**}$, or $M_*(\frak{G})^*$, then the space $\hbox{Der}_{p, q}(\frak{A})$ is trivial.
\end{corollary}

\section{\normalsize\bf $ (p,q)$-Jordan derivations of group algebras }

An additive mapping $\frak{D}: A\rightarrow A $ is called a $ \textit{$(p,q)$-Jordan derivation} $ if
\begin{eqnarray}\label{r2}
(p+q) \frak{D}(a^{2}) = 2p \frak{D}(a)a + 2q a \frak{D}(a) 
\end{eqnarray}
for all $ a\in A $. Note that if $p=0$ and $q>0$ (respectively, $p>0$ and $q=0$), then $\frak{D}$ is called a \emph{Jordan left} (respectively, \emph{right}) \emph{derivation}. We denote by $\hbox{Der}_{p, q}^J(A)$, $\hbox{Der}_{l}^J(A)$ and $\hbox{Der}_{r}^J(A)$ the set of all $(p, q)$-Jordan derivations, Jordan left derivations and Jordan right derivations of $A$, respectively. Vukman $\cite{16}$ introduced the notion of $(p,q)$-Jordan derivations and showed that  
$$
\hbox{Der}_{p, q}^J(R)\subseteq\hbox{Der}(R),
$$
where $p$ and $q$ are two distinct positive integers and $R$ is a prime ring with a suitable characteristic. Kosi-Ulbl and Vukman $\cite{9}$ established this result for semiprime rings; see also \cite{sa, bdf, fv}.  In the next result we investigate $\hbox{Der}_{p, q}^J(L_0^\infty(\frak{G})^*)$. To this end, let us recall that the right annihilator of $L_0^\infty(\frak{G})^*$, denoted by $\hbox{ran}(L_0^\infty(\frak{G})^*)$, is the set of all $\gamma\in L_0^\infty(\frak{G})^*$ such that $F\diamond\gamma=0$ for all $F\in L_0^\infty(\frak{G})^*$. 

\begin{theorem}\label{habb} Let $p$ and $q$ be distinct positive integers. Then
$$
\hbox{Der}_{p, q}^J(L_0^\infty(\frak{G})^*)\subseteq \hbox{Der}_{l}^J(L_0^\infty(\frak{G})^*).
$$
Furthermore, the range of $\frak{D}\in \hbox{Der}_{p, q}^J(L_0^\infty(\frak{G})^*)$ is contained in $\emph{ran}(L_0^\infty(\frak{G})^*)\subseteq\emph{rad}(L_0^\infty(\frak{G})^*)$.
\end{theorem}
{\it Proof.} Let $\frak{D}\in \hbox{Der}_{p, q}^J(L_0^\infty(\frak{G})^*)$. Then for every $F\in L_0^\infty(\frak{G})^*$, we have
\begin{eqnarray}\label{hab}
(p+q)\frak{D}(F^2)=2p\frak{D}(F)\diamond F+2qF\diamond \frak{D}(F).
\end{eqnarray}
Thus 
$$
(p+q)\frak{D}(E)=2p\frak{D}(E)+2qE\diamond \frak{D}(E),
$$
where $E$ is a right identity for $L_0^\infty(\frak{G})^*$. This shows that 
$$
-pE\diamond \frak{D}(E)=2qE\diamond \frak{D}(E)-qE\diamond \frak{D}(E)=qE\diamond \frak{D}(E).
$$
Hence $E\diamond \frak{D}(E)=0$. So $\frak{D}(E)\in\hbox{ran}(L_0^\infty(\frak{G})^*)$. The  linearization (\ref{hab}) gives
\begin{eqnarray}\label{hab1}
(p+q)\frak{D}(F_1\diamond F_2+F_2\diamond F_1)&=&2p\frak{D}(F_1)\diamond F_2+2qF_1\diamond \frak{D}(F_2)\\
&+&2p\frak{D}(F_2)\diamond F_1+2qF_2\diamond \frak{D}(F_1)\nonumber
\end{eqnarray}
for all $F_1, F_2\in L_0^\infty(\frak{G})^*$. Set $F_1=E$ and $F_2=\gamma\in\hbox{ran}(L_0^\infty(\frak{G})^*)$ in (\ref{hab1}). Then
\begin{eqnarray*}
(p+q)\frak{D}(\gamma)&=&2qE\diamond \frak{D}(\gamma)+2p\frak{D}(\gamma)+2q\gamma\diamond \frak{D}(E)\\
&=&2qE\diamond \frak{D}(\gamma)+2p\frak{D}(\gamma).
\end{eqnarray*}
It follows that
$$
(q-p)\frak{D}(\gamma)=2qE\diamond \frak{D}(\gamma).
$$
Consequently,  $(q-p)E\diamond \frak{D}(\gamma)=2qE\diamond \frak{D}(\gamma)$ which implies that $\frak{D}(\gamma)\in\hbox{ran}(L_0^\infty(\frak{G})^*)$. We define $\bar{\frak{D}}: L_0^\infty(\frak{G})^*/\hbox{ran}(L_0^\infty(\frak{G})^*)\rightarrow L_0^\infty(\frak{G})^*/\hbox{ran}(L_0^\infty(\frak{G})^*)$ by $$\bar{\frak{D}}(F+\hbox{ran}(L_0^\infty(\frak{G})^*))=\frak{D}(F)+\hbox{ran}(L_0^\infty(\frak{G})^*).$$ 
It is easy to show that $\bar{\frak{D}}\in \hbox{Der}_{p, q}^J(L_0^\infty(\frak{G})^*)$. On the other hand, in view of Theorem 2.11 of \cite{lp}, $L_0^\infty(\frak{G})^*/\hbox{ran}(L_0^\infty(\frak{G})^*)$ and $M(\frak{G})$ are isometrically isomorphic. Hence $L_0^\infty(\frak{G})^*/\hbox{ran}(L_0^\infty(\frak{G})^*)$ is semisimple. By Theorem 2.5 of \cite{6}, $\bar{\frak{D}}$ is zero. That is, $\frak{D}$ maps $L_0^\infty(\frak{G})^*$ into $\hbox{ran}(L_0^\infty(\frak{G})^*)$. Hence 
\begin{eqnarray}\label{hab2}
(p+q)\frak{D}(F^2)=2p\frak{D}(F)\diamond F
\end{eqnarray}
for all $F\in L_0^\infty(\frak{G})^*$. Note that if $p=q$, then $\frak{D}\in \hbox{Der}_{l}^J(L_0^\infty(\frak{G})^*)$. Hence we may assume that $p\neq q$. Take $F=E$ in (\ref{hab2}). Then $\frak{D}(E)=0$. From (\ref{hab1}) and the fact that $\frak{D}(L_0^\infty(\frak{G})^*)\subseteq\hbox{ran}(L_0^\infty(\frak{G})^*)$ we see that 
\begin{eqnarray*}
(p+q)\frak{D}(E\diamond F+F)&=&2p\frak{D}(E)\diamond F+2qE\diamond \frak{D}(F)\\
&+&2p\frak{D}(F)\diamond E+2qF\diamond \frak{D}(E)\\
&=&2qE\diamond \frak{D}(F)+2p\frak{D}(F)\\
&=&2p\frak{D}(F)
\end{eqnarray*}
for all $F\in L_0^\infty(\frak{G})^*$. Thus 
\begin{eqnarray}\label{hab3}
\frak{D}(F)=\frac{p+q}{p-q}\frak{D}(E\diamond F).
\end{eqnarray}
The substitution $F-E\diamond F$ for $F$ in (\ref{hab3}) gives $\frak{D}(F)=\frak{D}(E\diamond F)$. This together with (\ref{hab3}) shows that $q=0$. Now from (\ref{hab2}) we conclude that $\frak{D}\in \hbox{Der}_{l}^J(L_0^\infty(\frak{G})^*)$.$\hfill\square$

\begin{corollary}\label{ali} Let  $\frak{G}$ be a locally compact abelian group. Then $\hbox{Der}_{p, q}^J(L_0^\infty(\frak{G})^*)=\{0\}$.
\end{corollary}
{\it Proof.} Let $\frak{D}\in \hbox{Der}_{p, q}^J(L_0^\infty(\frak{G})^*)$. If $p=0$, then $\frak{D}\in \hbox{Der}_{l}^J(L_0^\infty(\frak{G})^*)$ and so $\frak{D}=0$; see Theorem 4.1 of \cite{am}. If $q=0$, then $\frak{D}\in \hbox{Der}_{r}^J(L_0^\infty(\frak{G})^*)$ and by Theorem 3.1 of \cite{am}, $\frak{D}$ is a right derivation. Since each Jordan right derivation is a $(0, q)$-derivation, it follows from Corollary \ref{am} that $\frak{D}=0$. In the case where, $p$ and $q$ are distinct positive integers, Theorem \ref{habb} implies that $\frak{D}=0$.$\hfill\square$\\

Let $R$ be a prime ring with $\hbox{char}(R) \neq 2$. By Posner's first theorem \cite{p}, if $D_1, D_2, D_1D_2 \in \hbox{Der}(R)$, then either $D_1 = 0$ or $D_2 = 0$. Creedon \cite{c} studied a version of Posner's first theorem for algebras and showed that if the conditions of Posner's theorem hold, then
$$
D_1D_2(A)\subseteq \hbox{rad}(A).
$$
The first author together with Ahmadi \cite{am} studied $\hbox{Der}_{r}^J(L_0^\infty(\frak{G})^*)$ in the case where $\frak{G}$ is abelian, and posed the following open problem.

\begin{problem} {\rm Does Posner's first theorem remain valid for the elements of $\hbox{Der}_{r}^J(L_0^\infty(\frak{G})^*)$?}
\end{problem}

As a consequence of Corollary \ref{ali}, we provide an affirmative answer to this problem.

\begin{corollary}
Let $\frak{G}$ be a locally compact abelian group. Then
$$
\hbox{Der}_{r}^J(L_0^\infty(\frak{G})^*) = \{0\}.
$$
In particular, the product of any two elements $D_1, D_2 \in \hbox{Der}_{r}^J(L_0^\infty(\frak{G})^*)$ is identically zero.
\end{corollary}

\end{document}